\documentclass[12pt]{amsart}
\usepackage[cp1250]{inputenc}
\usepackage[T1]{fontenc}
\usepackage{amscd}

\theoremstyle{plain}
\newtheorem{thm}{Theorem}[section]
\newtheorem{lem}[thm]{Lemma}
\newtheorem{prop}[thm]{Proposition}
\newtheorem{cor}[thm]{Corollary}

\theoremstyle{definition}
\newtheorem{dff}[thm]{Definition}
\newtheorem{rem}[thm]{Remark}
\newtheorem{exa}[thm]{Example}

\numberwithin{equation}{section}
\def\a{\alpha}

\def\B{\mathcal{B}}

\def\D{\mathcal{D}}
\def\E{\mathcal{E}}

\def\H{\mathcal{H}}

\def\phi{\varphi}
\def\p{\partial}
\def\pp{\phi}

\def\ci{\circ}
\def\rz{\mathbb{R}}
\def\R{\rz}

\def\wyr#1{\textit{#1}}

\def\s{\subset}

\def\t{\times}
\def\r{\rightarrow}

\def\ld{,\ldots,}
\def\pp{\partial}

 \DeclareMathOperator{\diff}{Diff}
 \DeclareMathOperator{\cl}{cl}
  \DeclareMathOperator{\cld}{cld}
 \DeclareMathOperator{\Frag}{Frag}
 \DeclareMathOperator{\frag}{frag}
 \DeclareMathOperator{\fragd}{fragd}
 \DeclareMathOperator{\Fragd}{Fragd}
\DeclareMathOperator{\id}{id} \DeclareMathOperator{\fl}{Fl}
\DeclareMathOperator{\intt}{Int} \DeclareMathOperator{\supp}{supp}

\thanks{Partially supported by the Polish Ministry of Science and Higher Education and the
AGH grant n. 11.420.04}
\address{Faculty of Applied Mathematics, AGH University of Science and
\linebreak Technology, al. Mickiewicza 30, 30-059 Krak\'ow,
Poland} \email{tomasz@uci.agh.edu.pl}
\date{revised; June 7, 2010}

\title{Boundedness of certain automorphism groups of an open manifold}
\author{ Tomasz Rybicki}

\begin{document}

\maketitle

\begin{abstract}

It is shown that certain diffeomorphism or homeomorphism groups
with no restriction on support of an open manifold (being the
interior of a compact manifold) are bounded. It follows that these
groups are uniformly perfect. In order to characterize the
boundedness several conditions on automorphism groups of an open
manifold are introduced. In particular, it is shown that the
commutator length diameter of the automorphism group $\D(M)$ of a
portable manifold $M$ is estimated by 4.
\end{abstract}

\section{Introduction}

Let us recall that a group is called {\it bounded} if it is
bounded with respect to any bi-invariant metric. The purpose of
this paper is to show that some diffeomorphism or homeomorphism
groups with not necessarily compact support of an open manifold
are bounded. We will formulate some conditions which ensure the
boundedness of such groups. Throughout, to avoid complications in
terminology we will refer to the homeomorphisms as diffeomorphisms
(of class $C^0$).

In the sequel we will deal with a manifold $M$ being the interior
of a compact manifold $\bar M$. We adopt the following notation
similar to that in \cite{bib:MD}. Let $\p_i$, $i=1\ld k$, be the
family of all connected components of the boundary $\p$ of $\bar
M$. Let $K=\{1\ld k\}$. For any $J\s K$ and $r=0,1\ld\infty$ the
symbol $\diff^r_J(M)$ will stand for the totality of
$C^r$-diffeomorphisms which are equal to the identity on a
neighborhood of $\p_J:=\bigcup_{i\in J}\p_i$. Then
$\diff^r(M)=\diff^r_{\emptyset}(M)$ and
$\diff^r_c(M)=\diff^r_K(M)$, the group of compactly supported
$C^r$-diffeomorphisms of $M$. Next, $\D^r_J(M)$ denotes the group
of all elements of $\diff^r_J(M)$ that can be joined with the
identity by an isotopy in $\diff^r_J(M)$. In particular,
$\D^r(M)=\D^r_{\emptyset}(M)$ (resp. $\D^r_c(M)=\D^r_K(M)$) is the
(resp. compactly supported) identity component  of the group of
all $C^r$-diffeomorphisms of $M$.

 The
problem of the boundedness of a group of diffeomorphisms is
closely related to its uniform perfectness (c.f. \cite{bib:BIP},
Propositions 1.3 and 1.4). Recall that a group $G$ is called
\wyr{perfect} if it is equal to its own commutator subgroup
$[G,G]$. Next $G$ is said to be \wyr{uniformly perfect} if $G$ is
perfect and there exists a positive integer $N$ such that any
element of $G$ can be expressed as a product of at most $N$
commutators of elements of $G$.  For $g\in [G,G]$ the least $N$
such that $g$ is a product of $N$ commutators is called the
\wyr{commutator length} of $g$ and is denoted by $\cl_G(g)$.

If $\pp\bar M=\emptyset$ then it is well-known that $\D^r(M)$ is a
simple group, where $r=0,1\ld \infty$, except possibly
$r=\dim(M)+1$ (c.f. for $r=0$ \cite{bib:Mat71} with
\cite{bib:EdwKir71}, and for $r\geq 1$ \cite{bib:Thu74} and
\cite{bib:Mat74}). Another basic theorem was proved by D. McDuff
in \cite{bib:MD}.

\begin{thm} \cite{bib:MD}
Suppose that $\pp\bar M\neq\emptyset$. The groups $\D^r_J(M)$ are
perfect unless $J=K$ and $r=\dim(M)+1$. In particular, $\D^r(M)$
is a perfect group.
\end{thm}

In the sequel we will always assume that $\pp\bar M\neq\emptyset$.

In  2009 P. Schweitzer reconstructed in \cite{bib:Sch} his own
proof of the fact that the quotient
$\D^r(M)/\diff^r_{\{i\}}(M)\cap\D^r(M)$, $i\in K$, is a simple
group. This theorem, which is a "hardest work" (\cite{bib:MD}) in
the proof of Theorem 1.1, had been proved independently by W.
Ling, Schweitzer and McDuff more than thirty years ago, but the
proof has never been published.

Recently, basic results concerning the uniform perfectness and the
boundedness of diffeomorphism groups of many manifolds have been
proved  by D. Burago, S. Ivanov and L. Polterovich in
\cite{bib:BIP} and by T. Tsuboi in \cite{bib:Tsu2}.  In contrast
to the problem of perfectness and simplicity, these results depend
essentially on the topology of the underlying manifold. These
results generalize older ones, e.g. \cite{bib:An}.

Note that the problem of the uniform perfectness and the
boundedness is still valid for some nontransitive diffeomorphism
groups which are perfect but non-simple, e.g. for the
diffeomorphism group of manifold with boundary (\cite{bib:Ryb98}),
or  of a foliated manifold (\cite{bib:Ryb95}, \cite{bib:Tsu06},
\cite{bib:LR}). For the problem of perfectness and uniform
perfectness in the relative case of $\diff(M,N)$, where $N$ is a
proper submanifold of $M$, see \cite{bib:Ryb1998} and
\cite{bib:AF}.

 In section 2 we will show the
equivalence of many conditions describing automorphism groups of
an open manifold. In particular, as a consequence of these results
  we have
\begin{thm} Let $M$ be a manifold of class $C^r$ , $r=0,1\ld\infty$, as above and let $M$ be
portable or, more generally, let $M$ satisfy the
$(\sqcup)$-property (c.f. Def.2.14). Let $J\s K$ and
$r\neq\dim(M)+1$. Then the group  $\D^r_J(M)$ is bounded. In
particular, $\D^r(M)$ is bounded.
\end{thm}

Notice that the class of portable manifolds includes the euclidean
spaces $\R^n$, the manifolds of the form $M\t\R^n$, and the
three-dimensional handlebodies.

\begin{cor} Under the above assumptions,
the group $\D_J^r(M)$ is uniformly perfect.  The commutator length
diameter of $\D_J^r(M)$ is $\leq 4$.
\end{cor}

The first statement is  an immediate consequence  of Theorem 1.2,
the second will be proved in section 2. In section 2 we also
obtain another estimation on the diameter of $\cl_{\D^r_J(M)}$,
c.f. (2.2), for a much wider class of automorphism groups.

In our investigations a special role is played by  the properties
that an automorphism group is "factorizable" or that it is
"determined on compact subsets" (Definitions 2.1 and 2.5(2)). The
significance of these properties is illustrated by examples in
section 3.

{\bf Acknowledgments.} A correspondence with Paul Schweitzer and
his recent paper \cite{bib:Sch} were helpful when I was preparing
my paper. I would like to thank him for his kindness. I also
express my deep gratitude to the referee for pointing out some
unclear or wrong statements, especially  in a previous version of
Def. 2.1, and for other valuable  remarks.

\section{Boundedness and uniform perfectness of certain automorphism groups}

The notion of the conjugation-invariant norm is an indispensable
tool in studies of the boundedness of groups.  Let $G$ be a group.
A \wyr{conjugation-invariant norm} on $G$ is a  function
$\nu:G\r[0,\infty)$ which satisfies the following conditions. For
any $g,h\in G$ \begin{enumerate} \item $\nu(g)>0$ if and only if
$g\neq e$; \item $\nu(g^{-1})=\nu(g)$; \item
$\nu(gh)\leq\nu(g)+\nu(h)$; \item $\nu(hgh^{-1})=\nu(g)$.
\end{enumerate}

It is easily seen that $G$ is bounded if and only if any
conjugation-invariant norm on $G$ is bounded.

Suppose that $G$ is perfect. Then the commutator length $\cl_G$ is
a conjugation-invariant norm on $G$.

Recall that $M$ is the interior of a compact, connected manifold
$\bar M$ of class $C^r$, where $ r=0,1\ld\infty$, with non-empty
boundary $\p=\p_1\cup\ldots\cup\p_k$. By a \wyr{product (or
collar) neighborhood} of $\p$ we mean a closed subset
$P=\bigcup_{i=1}^k P^{(i)}$ of $M$, where $P^{(i)}=\p_i\t[0,1)$,
so that $P=\p\t[0,1)$. Here $\p\t[0,1]$ is embedded in $\bar M$,
and $\p\t\{1\}$ is identified with $\p$. %Denote $\breve
%P^{(i)}=\p_i\t(0,1)$ and $\breve P=\bigcup_{i=1}^k \breve
%P^{(i)}$.
In particular, $P^{(i)}$ are pairwise disjoint.

A \wyr{translation system} on the product manifold $N\t[0,\infty)$
(c.f. \cite{bib:Li1}, p.168) is a family $\{P_j\}_{j=1}^{\infty}$
of closed product neighborhoods of $N\t\{\infty\}$ such that
$P_{j+1}\s\intt P_j$ and $\bigcap_{j=1}^{\infty}P_j=\emptyset$. A
detailed description of the role played by translation systems on
product manifolds $N\t\R$ was given in Ling's paper
\cite{bib:Li1}.  By a translation system of $M$ we understand a
translation system on a product neighborhood of $\p_i$,
$P^{(i)}=\p_i\t[0,1)$, where   $i\in K$. By a {\it ball} we mean
an open ball with its closure compact and contained in a chart
domain.

Let $G$ be a subgroup of $\diff^r(M)$. For a subset $U\s M$ denote
by $G(U)$ the subgroup of all elements of $G$ which can be joined
with the identity by an isotopy in $G$ compactly  supported in
$U$. Recall that the support of a diffeomorphism $f$ on $M$,
denoted by $\supp(f)$, is the closure of the set $\{x\in M:
f(x)\neq x\}$, and the support of an isotopy $\{f_t\}_{t\in
[0,1]}$ is defined by $\supp(\{f_t\}):=\bigcup_t\supp(f_t)$.

We say that $g\in G$ \emph{meets the $i$-th end} of $M$ ($i\in K$)
if for any neighborhood $U$ of $\p_i$ we have
$U\cap\supp(g)\neq\emptyset$, i.e. if $g$ does not stabilize near
$\p_i$. Next, we say that $G$ meets the $i$-th end of $M$ if there
is $g\in G$ which does so. Denote $J_{G}=\{i\in K: G$\ meets the
$i$-th end of $M\}$, and\quad $\p_{G}=\bigcup_{i\in J_{G}}\p_i$.

\begin{dff}

 $G$ is called
\wyr{factorizable} if for any $g\in G$ there are a family of balls
$\{B_{\alpha}\}_{\alpha\in A}$, a product neighborhood
$P=\p\t[0,1)$, and  a family of diffeomorphisms $g_j\in G$ for
$j=0,1\ld N$ such that:

(1) $g=g_0 g_1\cdots g_N$ with $g_0\in G(P)$ and $g_j\in
G(B_{\alpha(j)})$, $j=1\ld N$.

Moreover, for any product neighborhood $P$ and for any $g\in G(P)$
there is a sequence of reals from (0,1) tending to 1
\[0<a_1<\bar a_1<\bar b_1<b_1<a_2<\ldots<a_n<\bar a_n<\bar
b_n<b_n<\ldots<1\] and $h\in G(P)$ such that

(2) $h=g$ on $\bigcup_{n=1}^{\infty} \p\t[\bar a_n,\bar b_n]$;

(3) $h=\id$ on $P^{(i)}$, $i\in K$, whenever $g=\id$ on $P^{(i)}$.

Put $D_n:=\p\t(a_n,b_n)$ and $D:=\bigcup_{n=1}^{\infty}D_n$. Then
we also assume that:

(4) $\supp(h)\s D$;

(5) for the resulting decomposition $h=h_1h_2\ldots$ with respect
to $D=\bigcup_{n=1}^{\infty}D_n$ we have $h_n\in G(D_n)$ for all
$n$.
\end{dff}

\begin{rem}
The reason for formulating Def. 2.1 is the absence of isotopy
extension theorems or fragmentation theorems for automorphism
groups of some geometric structures. Roughly speaking, $G$
satisfies Def. 2.1 if all its elements can be joined with id by an
isotopy in $G$ and appropriate versions of the above mentioned
theorems are available.
\end{rem}
Let $G$ be factorizable. Then for any $g\in G$ there are a family
of balls $\{B_{\alpha}\}$, a product neighborhood $P$ of $\p$, and
a decomposition \begin{equation} g=g_{01}\ldots g_{0k}g_1\ldots
g_N, \end{equation} where $g_{0i}\in G(P^{(i)})$, $i=1\ld k$, and
$g_j\in G(B_{\alpha_j})$, $j=1\ld N$. For $J\s K$ we define $G_J$
as the totality of $g\in G$ such that there is a decomposition
(2.1) with $g_{i0}=\id$ for all $i\in J$. We also put $G_c=G_K$.
Then it is easily checked the following

\begin{prop}
If $G$ is factorizable then so are $G_J$. Moreover, if $G$ is
factorizable then $G$ (resp. $G_J$) are contained in $\D^r(M)$
(resp. $\D^r_J(M)$). We also have   $J_{G_J}\s K\setminus J$.
\end{prop}

\begin{exa}
The group $\diff^r(\R^n)$ does not satisfy Def.2.1. The reason is
that in this case any $f\in\diff^r(\R^n)$ would be isotopic to id
due to 2.1(1) which is not true. On the other hand, any
$f\in\diff^r_c(\R^n)$ is isotopic  to the identity but the isotopy
need not be compactly supported. It follows that $\diff^r_c(\R^n)$
does not fulfil Def.2.1.(1). The exception is $r=0$, when the
Alexander trick is in use and any compactly supported
homeomorphism on $\R^n$ is isotopic to id by a compactly supported
isotopy. It follows that $\diff^0_c(\R^n)$ is factorizable in view
of \cite{bib:EdwKir71}.

Let $C=\R\t\mathbb S^1$ be the annulus. Then there is the twisting
number epimorphism $\diff^r_c(C)\r\mathbb Z$. It follows that
$\diff^r_c(C)$ is unbounded in view of Lemma 1.10 in
\cite{bib:BIP}. On the other hand, $\diff^r_c(C)$ is not
factorizable. (This part of example was suggested by the referee.)
\end{exa}
\begin{dff}
\begin{enumerate}

 \item $G$ is called {\it locally perfect} if there
exists a covering by balls $\B$ of $M$ such that for any ball
$B\in\B$ the subgroup $G(B)$ is perfect.

 \item $G$ is said to be \wyr{determined on compact subsets}
 if the following is satisfied. Let $f\in\D^r(M)$.
 If there are a sequence of relatively compact subsets  $U_1\s\overline U_1\s
 U_2\s\ldots\s U_n\s\overline{U}_n\s
 U_{n+1}\s\ldots$ with  $\bigcup U_n=M$ and  a sequence $\{g_n\}$, $n=1,2\ld$ of elements
 of $G$ such that $f|_{U_n}=g_n|_{U_n}$ for $n=1,2\ld$ then we
have $f\in G$.

\item We say that $G$ \wyr{admits translation systems} if for any
$i\in J_{G}$ and for any sequence $\{\lambda_n\}$, $n=0,1,\ldots$,
with $\lambda_n\in(0,1)$, tending increasingly to 1, there exists
a $C^r$-mapping $[0,\infty)\ni t\mapsto f_t\in G$ supported in the
interior of $P^{(i)}$, with $f_0=\id$, $f_j=(f_1)^j$ for
$j=2,3,\dots$, and such that for the translation system
$P_n=\p_i\t[\lambda_n,1)$ one has $f_1(P_n)=P_{n+1}$ for
$n=0,1,2,\ldots$.
\end{enumerate}
\end{dff}

The following is a version of Isotopy Extension Theorem (c.f.
\cite{bib:EdwKir71},  \cite{bib:hir}).
\begin{thm}
Let $f_t$ be an isotopy in $\D^r(M)$ and let $C\s M$ be a compact
set. Then for any open neighborhood $U$ of $\bigcup_{t\in
[0,1]}f_t(C)$ there is an isotopy $g_t$ in $\D^r(M)$ such that
$g_t=f_t$ on $C$ and $\supp(g_t)\s U$.
\end{thm}

\begin{prop}
The groups $\D^r_J(M)$ ($J\s K$) satisfy Definitions 2.1 and 2.5
with the possible exception for 2.5(1) in case $r=\dim(M)+1$. In
particular, if $G=\D^r(M)$ then $G_J=\D_J^r(M)$.
\end{prop}
\begin{proof}
The proof will be written for $\D^r(M)$ (for $\D^r_J(M)$ is the
same). For open $U\s M$ let $\D^r(U)$   denote the group of all
elements of $\D^r(M)$ that can be joined with the identity by an
isotopy supported in $U$.

Def. 2.1: Let $g\in\D^r(M)$ and let $\tilde g_t$ be an isotopy
joining $g$ with the identity. Take $\{B_{\alpha}\}$ and $P$ such
that for $U=\bigcup_{\alpha\in A}B_{\alpha}$ we have
$\overline{M\setminus P}\s U$. Choose a compact subset
$C=\p\t[\lambda,\mu]\s U$, where $0<\lambda<\mu<1$. Possibly
enlarging $U$, by  Theorem 2.6 there is an isotopy $f_t$ in
$\D^r(U)$ such that $f_t=\tilde g_t$ on $C$. Put
$h_t=f_t^{-1}\tilde g_t$. Then $h_t|_C=\id$ and $h_t=\tilde g_t$
off $U$. Set $h_t=\tilde h_tk_t$, where $\tilde
h_t\in\D^r(\intt(P))$, $\supp(k_t)\s U$, and $\supp(\tilde
h_t)\cap\supp(k_t)=\emptyset$. It follows that $\bar g_t=\tilde
h_t^{-1}\tilde g_t$ has support in $U$. Set $g_0=\tilde h_1$. By a
fragmentation property for isotopies (\cite{bib:Ban97},
\cite{bib:EdwKir71}) we get $\bar g_1=g_1\cdots g_N$ with $g_j\in
\D^r( B_{\alpha(j)})$ for $j=1\ld N$. Clearly $g=g_0 g_1\cdots
g_N$, hence (1).

To show (2)-(5) we apply Theorem 2.6. This enables us to define
recurrently $a_n<\bar a_n<\bar b_n<b_n$ and $h|_{[a_n,b_n]}$ for
$n=1,2,\ldots$ in such a way that the claim is fulfilled. In fact,
let $g_t$ be an isotopy joining $g$ with the identity and
supported in $\intt(P)$ and suppose we have defined $0<a_1<\ldots<
b_{n-1}$. It suffices to take $b_{n-1}<a_n<\bar a_n<\bar b_n<b_n$
in such a way that $\p\t[a_n,1)$ is disjoint with
$\bigcup_{t\in[0,1]}g_t^{-1}(\p\t[0,b_{n-1}])$ and use Theorem
2.6.

Def. 2.5: (1) It is a consequence of the fundamental results on
the simplicity of diffeomorphism groups (\cite{bib:Mat71},
 \cite{bib:Thu74}, \cite{bib:Mat74}, \cite{bib:EdwKir71}).
(2) is trivial.

(3) Let $\lambda_n>0$ be as above. Let
$\tau:[0,1]\r\diff^{\infty}([0,1))$ be an isotopy such that
$\tau_0=\id$  and $\tau_1(\lambda_n)=\lambda_{n+1}$ for
$n=0,1,2,\ldots$. Next, let $\tau:[1,2]\r\diff^{\infty}([0,1))$ be
an isotopy from $\tau_1$ to $\tau_2=(\tau_1)^2$. Continuing this
procedure with $\tau_j=(\tau_1)^j$, let
$\tau=\bigcup_{j=0}^{\infty}\tau|_{[j,j+1]}:[0,\infty)\r\diff^{\infty}([0,1))$,
where $\tau$ is smoothed on neighborhoods of  $j=1,2,\ldots$ if
necessary. Set $f_t=\id_{P^{(i)}}\t\tau_t$.
\end{proof}

Of course, if $g\in\diff_c^r(C)$, where $C=\R\t\mathbb S^1$, is a
diffeomorphism with nonzero twisting number, then $h$ determined
by Def. 2.1, (2)-(5), is no longer compactly supported and does
not extend to the boundary $\p$.

 The following fact
shows that the statement "$g_0\in G(P)$" in Def. 2.1(1) could be
replaced by a weaker one "$\supp(g_0)\s P$", provided $G$ fulfills
Def. 2.5(3).

\begin{prop}
Suppose that a group $G\s\diff^r(M)$  admits translation systems
(Def. 2.5(3)). If $g\in G$ with $\supp(g)\s P$, where $P$ is a
product neighborhood of $\p$, then $g\in \D^r(M)$.
\end{prop}
\begin{proof}
It suffices to show that any $g\in G$ with $\supp(g)\s\intt(P)$ is
isotopic to the identity. Take a sequence $\lambda_n>0$,
$n=0,1,\ldots$, tending increasingly to 1. We may arrange so that
$g|_{\p\t[0,\lambda_0]}=\id$.

Let an isotopy $f_t$ in $G$ be as in Def. 2.5(3). For $t\in(0,1]$
define
\[ g_t=f_{\frac{1-t}{t}}\ci g\ci f_{\frac{1-t}{t}}^{-1}.\]
Then $g_1=g$ and $g_t$ extends smoothly onto $[0,1]\t P$ so that
$g_0=\id$.
\end{proof}

It follows from the assumption on $M$ and Proposition 2.8 that
there is a compact set $C\s M$ such that if $f\in\diff^r(M)$ and
$f|C=\id$ then $f\in\D^r(M)$. Clearly, this is not true for all
open manifolds with finite number of ends.

\begin{lem}
If $G$ satisfies  Definitions 2.1 and 2.5, then any $g\in G(P)$,
where $P$ is a product neighborhood of $\p$, can be written as a
product of two commutators of elements of $G(P)$.

\end{lem}

\begin{proof}
We may assume that $g\in G(\intt( P))$. Choose as in Def. 2.1 a
sequence $0<a_1<\bar a_1<\bar b_1<b_1<a_2<\ldots<a_n<\bar a_n<\bar
b_n<b_n<\ldots<1$ and $h\in G(P)$ such that conditions (2)-(5) in
Def. 2.1 are fulfilled.%$h=g$ on $\bigcup_{n=1}^{\infty} \p\t[\bar
%a_n,\bar b_n]$ and $\supp(h)\s\bigcup_{n=1}^{\infty} \p\t( a_n,
%b_n)$.
Put $\bar h=h^{-1}g$, that is $g=h\bar h$. Then $\supp(\bar h)$ is
in $(0,\bar a_1)\cup\bigcup_{n=1}^{\infty}(\bar b_n,\bar
a_{n+1})$, and $\bar h=g$ on  $[0, a_1]\cup\bigcup_{n=1}^{\infty}[
b_n, a_{n+1}]$. It is easily observed that $\bar h$ also fulfills
(2)-(5) in Def. 2.1. It suffices to show that $h$ is a commutator
of elements in $G(\intt(P))$ (in the same way it is true for $\bar
h$).

Choose arbitrarily $\lambda_0\in (0,a_1)$ and
$\lambda_n\in(b_n,a_{n+1})$ for $n=1,2,\ldots$. In view of Def.
2.5(3) there exists an isotopy  $[0,\infty)\ni t\mapsto f_t\in G$
supported in $\p_{G}\t(0,1)$, such that $f_0=\id$ and
$f_j(P_n)=P_{n+j}$ for $j=1,2,\ldots$ and for $n=0,1,2,\ldots$,
where $P_n=\p_{G}\t[\lambda_n,1)$ for $n=0,1,\ldots$. Now define
$\tilde h\in G(\intt(P))$ as follows. Set $\tilde h=h$ on
$\p_{G}\t[0,\lambda_1)$, and $\tilde
h=h(f_1hf_1^{-1})\ldots(f_nhf_n^{-1})$ on
$\p_{G(M)}\t[0,\lambda_{n+1})$ for $n=1,2\ldots$. Here
$f_n=(f_1)^n$. Then $\tilde h|_{\p_G\t[0,\lambda_n)}$ is a
consistent family of functions, and $\tilde
h=\bigcup_{n=1}^{\infty} \tilde h|_{\p_G\t[0,\lambda_n)}$ is a
local diffeomorphism. It is easily checked that $\tilde h$ is a
bijection. Due to Def. 2.5(2) $\tilde h\in G(\intt (P))$.

By definition we have the equality $\tilde h=hf_1\tilde h
f_1^{-1}$. It follows that $h=\tilde h f_1 \tilde h^{-1}f_1^{-1}$,
as claimed.

\end{proof}

\begin{rem}  It is necessary to use a decomposition $g=h\bar h$ in
the above proof. In fact, we can proceed as above to define
$\tilde h$ directly from $g$ (instead from $h$), but we cannot
ensure that the resulting $\tilde h$ is surjective. On the other
hand, if we would try to define $\tilde h$ by using $\tilde
g_n=(f_ngf_n^{-1})\ldots(f_1gf_1^{-1})g$ on
$\p_{G(M)}\t[0,\lambda_{n+1})$ then the family $\tilde g_n$ is
inconsistent and we cannot glue-up $\tilde g_n$.

\end{rem}

Suppose that $G$ is factorizable. For $g\in G$ we define
  $\frag_{G}(g)$ as the smallest $N$ such that there are a family of balls
  $\{B_{\alpha}\}$, a product neighborhood $P$ and
  and a decomposition (2.1).  Then
   $\frag_{G}$ is a conjugation-invariant norm on $G$,
called the \wyr{fragmentation norm}. In fact, since $G\s\D^r(M)$,
any $g\in G$ does not change the ends of $M$ so that it takes (by
conjugation) any decomposition in the form (2.1) into another
decomposition in the same form.

Define $\fragd_{G}:=\sup_{g\in G}\frag_{G}(g)$, the diameter of
$G$ in $\frag_{G}$. Consequently, $\frag_{G}$ is bounded iff
$\fragd_{G}<\infty$.

If $g\in G_c$ then, as usual (\cite{bib:BIP}), we define
$\Frag_{G_c}(g)$ to be the smallest $N$ such that $g=g_1\ldots
g_N$ with $g_i\in G(B_i)$, $B_i$ being balls, $i=1\ld N$, and
$\Fragd_{G_c}:=\sup_{g\in G_c}\Frag_{G_c}(g)$.

\begin{prop} For all $J\s K$ we have
$$\fragd_{G}=\fragd_{G_J}=\fragd_{G_c}=\Fragd_{G_c}.$$
\end{prop}
\begin{proof}
If $g\in G_c$ then $\Frag_{G_c}(g)\geq\frag_{G_c}(g)$ as any
fragmentation of $g$ supported in balls is of the form (2.1). On
the other hand, if $g=g_0g_1\ldots g_{N'}$ with
$N'<N=\Frag_{G_c}(g)$ is of the form (2.1), then $g_0^{-1}g\in
G_c$ and $\Frag_{G_c}(g_0^{-1}g)\leq N'$. Hence
$\fragd_{G_c}=\Fragd_{G_c}$.

Next, let $g=g_0g_1\ldots g_N\in G$ and let
$N=\frag_{G}(g)=\fragd_{G}$.  Take $h=g_1\ldots g_N\in G_c$. Then
$\frag_{G_c}(h)\leq\Frag_{G_c}(h)\leq N$. But if
$\frag_{G_c}(h)<N$ then $\frag_{G}(g)<N$, a contradiction.
Therefore, $\frag_{G_c}(h)=\frag_{G}(g)$. It follows that
$\fragd_G=\fragd_{G_c}$, since for $g\in G_c$ we have trivially
$\frag_G(g)=\frag_{G_c}(g)$.  Analogous statements are true for
$G_J$ instead of $G_c$.
\end{proof}

For any perfect group $G$ denote by $\cld_G$ the commutator length
diameter of $G$, i.e. $\cld_G:=\sup_{g\in G}\cl_G(g)$. Then $G$ is
uniformly perfect iff $\cld_G<\infty$.

Summing-up the above facts we have the following  generalization
of Theorem 1.1.
\begin{thm}
Assume that $G\s\diff^r(M)$ satisfies  Definitions 2.1 and 2.5.
Then $G$ and $G_J$ ($J\s K$)  are perfect groups. Moreover, if
there is a positive integer $r$ such that $\cld_{G(B)}\leq r$ for
all balls $B$ and the fragmentation norm $\frag_{G}$ of $G$ is
bounded, then $G$  and $G_J$ are uniformly perfect and we have the
inequalities
\begin{equation}
\cld_{G}\leq r\fragd_{G}+2,\quad \cld_{G_J}\leq
r\fragd_{G}+2.\end{equation}
\end{thm}

Observe that Theorem 2.12 is not  true for
$\diff^r_J(M)\cap\D^r(M)$ whenever $J$ is nonempty. In fact, for
$J\s K$, one has
\[[\diff^r_J(M)\cap\D^r(M),\diff^r_J(M)\cap\D^r(M)]=\D^r_J(M),\]
 unless $J=K$ and
$r=\dim(M)+1$, thanks to McDuff \cite{bib:MD}.

We will need some algebraic tools which mimic classical tricks for
homeomorphism groups (see, e.g., \cite{bib:An}). A subgroup $H$ of
$G$ is called  \wyr{strongly m-displaceable} if there is $f\in G$
such that the subgroups $H$, $fHf^{-1}$\ld $f^mHf^{-m}$ pairwise
commute. Then we say that $f$ \wyr{m-displaces} $H$. Fix a
conjugation-invariant norm $\nu$ on $G$ and assume that $H\s G$ is
strongly $m$-displaceable. Then $e_m(H):=\inf\nu(f)$, where $f$
runs over the set of elements of $G$ that $m$-displaces $H$, is
called the \wyr{order m displacement energy} of $H$.

Now in view of Theorem 2.2 in \cite{bib:BIP} we have that given
$G$ equipped with a conjugation-invariant norm $\nu$ and given
$H\s G$, a strongly $m$-displaceable subgroup of $G$, for any
$h\in[H,H]$ with $\cl_H(h)=m$ one has
\begin{equation}
\nu(h)\leq 14e_m(H)
\end{equation}
and
\begin{equation}
\cl_G(h)\leq 2.
\end{equation}
If $m=1$, i.e. $h$ is a commutator of elements of $H$,  then
\begin{equation}
\nu(h)\leq 4e_1(H).
\end{equation}
In particular, if there exists $g\in G$ that $m$-displaces $H$ for
every $m\geq 1$ the inequality (2.3) yields for all $h\in [H,H]$
that \begin{equation} \nu(h)\leq 14\nu(g).
\end{equation}

 A group $G\s\D^r(M)$ is said to be {\it locally moving} if for any
ball $B\s M$ and any $x\in B$ there is $f\in G(B)$ with $ f(x)\neq
x$. Next, $G$ acts {\it transitively inclusively} on $M$ if  for
every balls $U$ and $V$ there is $f\in G$ with $f(U)\s V$. It is
clear that $\D^r(M)$ is locally moving and acts transitive
 inclusively on the basis of all balls.
\begin{prop}
Assume  that $G$ is locally moving and $G$ acts transitively
inclusively on $M$. Then 0 is not an accumulation point in the set
of values of any conjugation-invariant norm on $G$.
\end{prop}
\begin{proof}
First observe the following fact:

$(*)$ For any ball $U$ there are non-commuting $f_1,f_2\in G(U)$.

Indeed, take $x\in U$ and $f_1\in G(U)$ such that $f_1(x)\neq x$.
Next choose $V\s U$ with $x\in V$ and $f_1(x)\not\in V$. Let
$f_2\in G(V)$ such that $f_2(x)\neq x$. It follows that
$f_2(f_1(x))=f_1(x)\neq f_1(f_2(x))$.

Fix $U$ and $f_1,f_2$ as in $(*)$. Let $\nu$ be a
conjugation-invariant norm such that for any $\epsilon>0$ there is
$g\in G$ with $0<\nu(g)<\epsilon$. As $g\neq\id$ it follows the
existence of  a ball $B$ with $g(B)\cap B=\emptyset$. Since $G$
acts transitively inclusively  there is $h\in G$ with $h(U)\s B$.
Therefore $h^{-1}gh(U)\cap U=\emptyset$. It follows from (2.5)
that
\[ \nu([f_1,f_2])\leq 4\nu(h^{-1}gh)=4\nu(g)<4\epsilon,\]
a contradiction.
\end{proof}

According to the terminology in \cite{bib:BIP} a group $G$ is
called \wyr{meagre} if it is bounded and discrete. The latter
means that 0 is not an accumulation point in the set of values of
any conjugation-invariant norm on $G$.

\begin{dff} (c.f. \cite{bib:BIP})
 A smooth connected
open manifold $M$ (with $M=\intt(\bar M)$, where $\bar M$ is a
compact manifold) is called \wyr{portable}
  if it admits a complete vector field $X$ and a compact set $C\s M$,
 called a \wyr{core} of $M$, such that the following conditions are
 satisfied:
 \begin{enumerate}
  \item for any compact set $K\s M$ there is $t>0$ such that $\fl^X_t(K)\s C$ (here
  $\fl^X_t$ is the flow of $X$);
\item there is $f\in \D^r_c(M)$ such that $f(C)\cap C=\emptyset$.
\end{enumerate}
More generally, a connected open manifold $M$ of class $C^r$,
$r=0,1\ld\infty$, satisfies \wyr{$(\sqcup)$-property} if there are
disjoint open subsets $U$, $V$ of $M$ such that there is $f\in
\D^r_c(M)$ with the closure of $f(U\cup V)$ contained in $V$, and
such that  for every $g_1\ld g_l\in \D_c^r(M)$ there is $h\in
\D^r_c(M)$ satisfying
\begin{equation}
h(\bigcup_{i=1}^l\supp(g_i))\s U.\end{equation}

\end{dff}
It is immediate that any portable manifold satisfies the
$(\sqcup)$-property.

\begin{prop}
If $U,V$ are open disjoint subsets of $M$ such that there is
$f\in\D^r(M)$ with $\overline{f(U\cup V)}\s V$ then $f$
$m$-displaces $\D^r(U)$ for all $m\geq 1$.

\end{prop}
Indeed, this follows from the relation $f^m(U)\s
f^{m-1}(V)\setminus f^m(V)$ for all $m\geq 1$.

The class of portable manifolds comprises the euclidean spaces
$\R^n$, the manifolds of the form $M\t\R^n$, or the manifolds
admitting an exhausting Morse function with finite numbers of
critical points such that all the indices are less that
$\frac{1}{2}\dim M$. In particular, every three-dimensional
handlebody is a portable manifold.

It will be useful a more general notion concerning homeomorphism
groups rather than manifolds.

\begin{dff}
A  group $G\s\D^r(M)$ on a manifold $M$ being the interior of a
compact manifold is said to satisfy \wyr{$(\E)$-property} if the
following conditions hold:
\begin{enumerate}
\item  $G$ acts
 transitively inclusively  on $M$.
 \item There are a ball
$B$, an open subset $U\s M$ disjoint with $B$ and $f\in G_c$ such
that the closure of $f(U\cup B)$ is contained in $U$. \item If
$j\in J_{G}$ and $P$ is a product neighborhood of $\p$, then for
any sequence in (0,1), tending to 1, of the form
\[0<a_1<b_1<a_2<b_2<\ldots<a_n<b_n<\ldots<1\]  there are $f_1,f_2\in G(\intt(P^{(j)}))$
 such that for $i=1,2,\ldots$ one has
\[f_1(\p_j\t([a_{2i-1},b_{2i-1}]\cup[ a_{2i}, b_{2i}]))\s \p_j\t(
a_{2i}, b_{2i}),\]
\[f_2(\p_j\t([a_{2i},b_{2i}]\cup[ a_{2i+1}, b_{2i+1}]))\s \p_j\t(
a_{2i+1}, b_{2i+1}).\] Moreover, if we have another sequence
tending to 1
\[0<\tilde a_1<\tilde b_1<\tilde a_2<\tilde b_2<\ldots<\tilde
a_n<\tilde b_n<\ldots<1\]   then there is an element of $G(\intt(
P^{(j)}))$ of the form $\id\t\phi$, where $\phi:[0,1)\r[0,1)$ is a
diffeomorphism, with $\phi(a_i)=\tilde a_i$ and $\phi(b_i)=\tilde
b_i$ for $i=1,2,\ldots$.
\end{enumerate}
\end{dff}

\begin{prop}
The groups  $\D^r_J(M)$ ($J\s K$) satisfy the $(\E)$-property
(Def. 2.16).
\end{prop}
The proof is obvious.

\begin{thm}  Let $J\s K$. Suppose that $G\s\diff^r(M)$ satisfies Definitions 2.1, 2.5 and
2.16.  Then  the following conditions are equivalent:
\begin{enumerate}
\item the norm $\frag_{G}$ is bounded; \item $G$ is bounded;
  \item $G_J$ is bounded.
  \end{enumerate}
 If $G$ is also locally moving  then the above  conditions  are equivalent to the
meagerness of  $G$, or of $G_J$.
\end{thm}

\begin{proof}
As $\frag_{G}$ is a conjugation-invariant norm the implication
(2)$\Rightarrow$(1) is trivial. To show (1)$\Rightarrow$(2)
suppose that $\frag_{G}$ is bounded. For any $g\in G$ we have a
decomposition $g=g_0g_1\cdots g_N$ specified in Def. 2.1 with $N$
bounded. In particular, there is  a family of balls
$\{B_{\alpha}\}_{\alpha\in A}$ and a product neighborhood $P$ of
$\p$ such that  $g_0\in G(P)$ and $g_j\in G(B_{\alpha(j)})$ for
$j=1\ld N$.

Moreover, for  $g_0\in G(P)$  there is a sequence, converging to
1, of the form
\[0<a_1<\bar a_1<\bar b_1<b_1<a_2<\ldots<a_n<\bar a_n<\bar
b_n<b_n<\ldots<1\] and $h_1,h_2\in G(P)$ such that $h_1=g$ on
$\bigcup_{n=1}^{\infty} \p\t[\bar a_{2n-1},\bar b_{2n-1}]$,
$\supp(h_1)\s U_1:=\bigcup_{n=1}^{\infty} \p\t( a_{2n-1},
b_{2n-1})$, $h_2=g$ on $\bigcup_{n=1}^{\infty} \p\t[\bar
a_{2n},\bar b_{2n}]$, and $\supp(h_2)\s
U_2:=\bigcup_{n=1}^{\infty} \p\t( a_{2n}, b_{2n})$. Now, applying
the reasoning from the proof of Lemma 2.9 for $h=h_1h_2$, it can
be checked that $g_0$ can be written as $g_0=h_1h_2h_3h_4$, where
$h_3=g$ on $\bigcup_{n=1}^{\infty} \p\t[ b_{2n-1}, a_{2n}]$,
$\supp(h_3)\s U_3:=\bigcup_{n=1}^{\infty} \p\t( \bar b_{2n-1},
\bar a_{2n})$, $h_4=g$ on $\bigcup_{n=0}^{\infty} \p\t[ b_{2n},
a_{2n+1}]$, and $\supp(h_4)\s U_4:=\bigcup_{n=0}^{\infty} \p\t(
\bar b_{2n}, \bar a_{2n+1})$, for some $0< \bar b_0< b_0< a_1$.
 Furthermore, $h_j$ satisfy the corresponding conditions (2)-(5) in
 Def. 2.1.

It follows  from Def. 2.16(2) and Proposition 2.15 (applied to
$G$) the existence of a ball $B$ and $f\in G_c$ such that
$\{f^m(B)\}_{m=0}^{\infty}$ is a pairwise disjoint family (here
$f^0=\id$). In view of Def. 2.16(1) there are $h_{\a(j)}\in G_c$
such that $h_{\a(j)}(B_{\a(j)})\s B$ for $j=1\ld N$. It follows
that $\{h_{\a(j)}^{-1}f^mh_{\a(j)}(B_{\a(j)})\}_{m=0}^{\infty}$ is
a pairwise disjoint family. Consequently, $h_{\a(j)}^{-1}f
h_{\a(j)}$ $m$-displaces $G(B_{\a(j)})$ for all  $m\geq 1$ and for
$j=1\ld N$.

 Likewise, in view of
Def. 2.16(3) and Proposition 2.15  we have the existence of
$f_{j}\in G$ with $\supp(f_j)\s\p_{G}\t(0,1)$ such that $G(U_j)$
is $m$-displaceable by $f_{j}$ for $j=1,2,3,4$ and for all $m\geq
1$.

Let $\nu$ be a conjugation-invariant norm on $G$. In view of (2.6)
and the invariance of $\nu$ we have
\begin{align*}
\nu(g)&\leq\nu(h_1)+\cdots+\nu(h_4)+\nu(g_1)+\cdots+\nu(g_N)\\
&\leq 14(\nu(f_1)+\cdots\nu(f_4)+N\nu(f)).\end{align*}

Although the sets $U_1\ld U_4$ depend on $g_0$ (and on $g$),
thanks to the second assertion of Def. 2.16(3) and the invariance
of $\nu$, the norms  $\nu(f_j)$ are independent of $g$. It follows
 that $\nu(g)$ is bounded, as required.

In view of Proposition 2.11 we have that (1) is equivalent to (3)
in the same way.
 The second assertion is a consequence
of Proposition 2.13.
\end{proof}

\begin{cor}
Let $\D(M)=\D^r(M)$, $r=0,1\ld \infty, r\neq \dim(M)+1$.
 The following conditions are equivalent:
\begin{enumerate}
\item the norm $\frag_{\D(M)}$ is bounded;  \item $\D(M)$ is
bounded;  \item $\D_J(M)$ is bounded; \item $\D_c(M)$ is bounded;
\item $\D(M)$ is meagre;\item $\D_J(M)$ is meagre; \item $\D_c(M)$
is meagre.
\end{enumerate}
\end{cor}
 It is a consequence of Propositions 2.7, 2.13 and 2.17, and of
Theorem 2.18.

\medskip

 \wyr{Proof of Theorem 1.2 and Corollary 1.3}.
Let $M$ satisfy the ($\sqcup$)-property. In view of Propositions
2.7 and 2.17, the group $\D^r(M)$ fulfills Definitions 2.1, 2.5
and 2.16 if $r\neq \dim(M)+1$. By (2.7), Proposition 2.15 and
(2.6) it is easy to check that $\D^r_c(M)$ is bounded. Therefore,
by Corollary 2.19 the groups $\D^r(M)$ and $\D^r_J(M)$ are
bounded. In particular, these groups are uniformly perfect.  The
inequality $\cld_{\D^r_J(M)}\leq 4$ follows from Proposition 2.7,
Lemma 2.9, and the fact that $\cld_{\D^r_c(M)}\leq 2$ (Theorem
1.18 with Remark 3.2 in \cite{bib:BIP}).

\section{Examples}

First recall the following basic fact, c.f. \cite{bib:BIP}.
\begin{prop}

Let $G$ be any group. If $H_1(G)$ is infinite then $G$ is
unbounded.
\end{prop}

The first  example reveals the significance of the property that
$G(M)$ is "determined on compact subsets" (Def. 2.5(2)).
\begin{exa}
 Let $N$ (resp. $S$) be the north (resp. south) pole of $\mathbb
 S^m$, and let $p_N:\mathbb S^m\setminus\{N\}\r\R^m$ (resp. $p_S:\mathbb
 S^m\setminus\{S\}\r\R^m$) be the corresponding stereographic
 projection. By $\D^{\infty}_c(\R^m,0)$ we denote the identity
 component of compactly supported $C^{\infty}$-diffeomorphisms of $\R^m$ fixing the
 origin $0\in\R^m$. Define $G(\mathbb
 S^m)=p_N^{-1}\D^{\infty}_c(\R^m,0)p_N$. Then $G(\mathbb
 S^m)\s\D^{\infty}(\mathbb S^m)$ and any element of $G(\mathbb
 S^m)$ preserves $S$. Next,  define
 $G(\R^m)=p_SG(\mathbb
 S^m)p_S^{-1}$.
 It is easily checked that $G(\R^m)\cong\D^{\infty}_c(\R^m,0)$ and
 $G(\R^m)$ fulfils Def. 2.1 and  (1), (3)  in Def. 2.5. Notice that  Def.
 2.5(2)
 does not hold. In fact, it suffices to take any
 $g\in\D^{\infty}(\mathbb S^m\setminus S)$ which is not
 smoothly extendable on $\mathbb S^m$, and $f=p_Sgp_S^{-1}$. Then
 $f\notin G(\R^m)$ but $f$ fulfills the assumption of Def. 2.5(2).
  Note that $G(\R^m)$ does not satisfy (the last statement) of Def.
  2.16(3).

 On the other hand, since
 $H_1(G(\R^m))=H_1(\D^{\infty}_c(\R^m,0))=\R$  by a theorem of K. Fukui in \cite{bib:Fuk80}, the
 group $G(\R^m)$ is unbounded by
 Proposition 3.1.

 Similarly, let $\D^{\infty}_{l,c}(\R^m,0)$ be the identity
 component of compactly supported $C^{\infty}$-diffeomorphisms of $\R^m$ fixing the
 origin and $l$-tangent to the identity at the origin ($l\geq 1$).
 In the same manner as above we define a group $G_l(\R^m)$ which
 satisfies Def. 2.1 and Def. 2.5, (1) and (3), but not Def. 2.5(2) and not Def. 2.16(3). In view of
 \cite{bib:Fuk80} we have
 $H_1(G_l(\R^m))=H_1(\D^{\infty}_{l,c}(\R^m,0))=\R^{l+1}$ so
 that the group $G_l(\R^n)$ is unbounded due to Proposition 3.1.
\end{exa}

A conjugation-invariant norm $\nu$ on a group $G$ is  \wyr{stably
bounded} if the limit $\lim_{n\r\infty}\frac{\nu(g^n)}{n}=0$ for
any $g\in G$. It is clear that any bounded $\nu$ is also stably
bounded. Next, a map $\phi:G\r\R$ is called a \wyr{quasi-morphism}
if there is $K>0$ such that $|\phi(gh)-\phi(g)-\phi(h)|\leq K$ for
any $g,h\in G$. A deep theorem of C. Bavard \cite{bib:bav} states
that the commutator length $\cl_G$ is stably unbounded if and only
if there exists a non-trivial (i.e. not being a morphism)
homogeneous quasi-morphism on $G$.
\begin{exa}
Let $\H(\R)$ be the group of all homeomorphisms $h$ of $\R$ which
verifies $h(x+1)=h(x)+1$ for all $x\in\R$. In view of
\cite{bib:EH} $h\in\H(\R)$ is a product of $p$ commutators if and
only if \begin{equation}
\inf_{x\in\R}(h(x)-x)<2p-1\quad\hbox{and}\quad
\sup_{x\in\R}(h(x)-x)>1-2p.\end{equation}
 The group $\H(\R)$
admits a quasi-morphism $\tau$, called the \wyr{translation
number}, given for $h\in\H(\R)$  by the formula
\[\tau(h)=\lim_{n\r\infty}\frac{h^n(x)-x}{n},\]
which is independent of $x\in\R$ thanks to (3.1), c.f.
\cite{bib:bav}. By the theorem of Bavard it follows that the
commutator length is stably unbounded on the commutator subgroup
$[\H(\R),\H(\R)]$ so that the commutator length on
$[\H(\R),\H(\R)]$ is unbounded. In view of Proposition 1.4 in
\cite{bib:BIP}, the group $\H(\R)$ is itself unbounded.

This is encoded in the fact that  Def. 2.1 is not fulfilled by
$\H(\R)$, while  Def. 2.5 holds. Indeed, Def. 2.5(1) is satisfied
due to \cite{bib:Mat71} and Def. 2.5(2) is trivial. To show Def.
2.5(3) we identify $[0,1)$ with $[0,\infty)$ and define
$[0,\infty)\ni t\mapsto f_t\in\H(\R)$ to be a continuous curve of
translations.

On the other hand, according to classical results of A. Denjoy
\cite{bib:Den} $\tau$ restricted to the subgroup in $\H(\R)$ of
strictly increasing smooth diffeomorphisms is trivial so that the
theorem of Bavard does not apply in this case. Notice that
$\D^{\infty}(\mathbb S)$, which is isomorphic to the subgroup of
strictly increasing smooth diffeomorphisms in $\H(\R)$, is bounded
due to Theorem 1.11 in \cite{bib:BIP}.
\end{exa}

Considering some subgroups of the symplectomorphism group or of
the volume preserving diffeomorphism group  may be a source of
examples of automorphism groups of open manifold that fulfill
 Def. 2.1 and Def. 2.5, (1) and (2), but not Def. 2.5(3). For instance, let
$G(\R^{2m})$ be the kernel of Calabi homomorphism of $\R^{2m}$
equipped with the standard symplectic form. It is well-known that
$G(\R^{2m})$ is a simple group by a classical result of A.
Banyaga, c.f. \cite{bib:Ban97}. Recently, D. Kotschick in
\cite{bib:Kot} proved that this group is also stably bounded.

Observe that $G(\R^{2m})$ is a compactly supported group. It is
not difficult to extend $G(\R^{2m})$ to a symplectomorphism group
$\hat G(\R^{2m})$ with no restriction on support by making use of
 Def. 2.1. Then $\hat G(\R^{2m})$ satisfies Def. 2.1 and Def. 2.5, (1) and (2), but not Def. 2.5(3).
The reason is that the symplectic volume is a symplectic
invariant. Consequently, Kotschick's criterion for stable
boundedness of the commutator length (Theorem 2.3 in
\cite{bib:Kot}) cannot be extended to groups with no restriction
on support.

\end{document}